%% file: sb_main.tex
\documentclass[12pt,english]{article}
\usepackage{fullpage}
\usepackage{babel}
\usepackage{exscale,xspace}
\usepackage{amsfonts,amssymb}

\usepackage[
  pdfauthor={Bernd C. Kellner},
  pdfkeywords={Bernoulli numbers, irregular pairs of higher order, index of irregularity,
    Kummer congruences, Adams theorem, Riemann zeta function, p-adic zeta function},
  pdftitle={The structure of Bernoulli numbers},
  colorlinks,bookmarksnumbered]{hyperref}

\include{sb_def}

\begin{document}

\title{The structure of Bernoulli numbers}
\author{Bernd C. Kellner}
\date{}
\maketitle

\abstract{We conjecture that the structure of Bernoulli numbers
can be explicitly given in the closed form
$$ B_n = (-1)^{\frac{n}{2}-1} \prod_{p-1 \notdiv n} |n|_p^{-1}
\prod\limits_{(p,l)\in\Psi^{\rm irr}_1 \atop n \equiv l \,\mods{p-1}} \!\!
|p \, (\chi_{(p,l)} - {\textstyle \frac{n-l}{p-1}})|_p^{-1} \,\,
\prod\limits_{p-1 \pdiv n} p^{-1} $$
where the $\chi_{(p,l)}$ are zeros of certain $p$-adic zeta functions
and $\Psi^{\rm irr}_1$ is the set of irregular pairs.
The more complicated but improbable case where the conjecture does not hold
is also handled; we obtain an unconditional structural formula for Bernoulli numbers.
Finally, applications are given which are related to classical results.}
\smallskip

\textbf{Keywords:} Bernoulli number, Kummer congruences, irregular prime,
irregular pair of higher order, Riemann zeta function, $p$-adic zeta function
\smallskip

\textbf{Mathematics Subject Classification 2000:} 11B68

\parindent 0cm

\input{sb_intro}
\input{sb_prelim}
\input{sb_delta}

\input{sb_sing}

\input{sb_appl}

\bigskip
Bernd C. Kellner \\
{\small
address: Reitstallstr. 7, 37073 G\"ottingen, Germany \\
email: bk@bernoulli.org}

\bibliographystyle{alpha}
\bibliography{sb_bib}

\end{document}

%% file: sb_def.tex
\newcommand{\binom}[2]{ {#1 \choose #2} }
\newcommand{\pdiv}{\mid}
\newcommand{\pdiveq}{\mid\mid}
\newcommand{\notdiv}{\nmid}
\newcommand{\sep}{:\,}

\newcommand{\mod}[1]{({\rm mod\ } #1)}
\newcommand{\mods}[1]{({\rm mod\,} #1)}
\newcommand{\qed}{\parbox{0cm}{}\hspace*{\fill} $\Box$}

\def\ord{\mathop{\rm ord}\nolimits}
\def\sign{\mathop{\rm sign}\nolimits}
\def\denom{\mathop{\rm denom}\nolimits}
\def\real{\mathop{\rm Re}\nolimits}

\newcommand{\CC}{\mathbb{C}}
\newcommand{\QQ}{\mathbb{Q}}
\newcommand{\ZZ}{\mathbb{Z}}
\newcommand{\NN}{\mathbb{N}}

\newcommand{\PP}{\mathbb{P}}

\newcommand{\refeqn}[1]{(\ref{#1})}

\newcommand{\IRR}{\Psi^{\rm irr}}
\newcommand{\IRP}{\widehat{\Psi}^{\rm irr}}
\newcommand{\eulerphi}{\varphi}
\newcommand{\Bnum}{\mathit{\Lambda}}
\newcommand{\Bden}{V}

\renewcommand{\theequation}{\arabic{section}.\arabic{equation}}
\newcommand{\settheequation}[1]{\renewcommand{\theequation}{#1}}
\newcommand{\restoretheequation}{\renewcommand{\theequation}%
{\arabic{section}.\arabic{equation}}\addtocounter{equation}{-1}}

\newtheorem{prop}{Proposition}[section]
\newtheorem{theorem}[prop]{Theorem}
\newtheorem{corl}[prop]{Corollary}

\newtheorem{conj}[prop]{Conjecture}

\newtheorem{thdefin}[prop]{Definition}
\newtheorem{thremark}[prop]{Remark}
\newtheorem{thtable}[prop]{Table}

\newenvironment{defin}{\begin{thdefin}\rm}{\end{thdefin}}
\newenvironment{remark}{\begin{thremark}\rm}{\end{thremark}}

\newenvironment{proof}{\medskip\noindent%
\textsc{Proof.}}{\qed\medskip}

{\end{equation}\restoretheequation}

%% file: sb_intro.tex
\section{Introduction}

The classical Bernoulli numbers $B_n$ are defined by the power series
\[
   \frac{z}{e^z-1} = \sum_{n=0}^\infty B_n \frac{z^n}{n!} \,,
     \qquad |z| < 2 \pi \,,
\]
where all numbers $B_n$ are zero with odd index $n > 1$.
The even-indexed rational numbers $B_n$ alternate in sign.
First values are given by $B_0 = 1$, $B_1 = -\frac{1}{2}$,
$B_2 = \frac{1}{6}$, $B_4 = -\frac{1}{30}$.
Although the first numbers are small with $|B_n| < 1$ for $n=2,4,\ldots,12$,
these numbers grow very rapidly with $|B_{n}| \rightarrow \infty$ for even
$n \rightarrow \infty$.
\medskip

For now, let $n$ be an even positive integer.
An elementary property of Bernoulli numbers is the following discovered
independently by T.~Clausen \cite{clausen40bern} and K.\,G.\,C. von Staudt \cite{staudt40bern}
in 1840. The structure of the denominator of $B_n$ is given by
\begin{equation} \label{eqn-clausen-von-staudt}
   B_n + \sum_{p-1 \pdiv n} \frac{1}{p} \in \ZZ \qquad \mbox{and} \qquad
   \denom(B_n) = \prod_{p-1 \pdiv n} p \,.
\end{equation}

As a result, often attributed to J.~C.~Adams, see \cite[Prop.~15.2.4, p.~238]{ireland90classical},
$B_n/n$ is a $p$-integer for all primes $p$ with $p-1 \notdiv n$. Therefore, one has a
trivial divisor of $B_n$ which is cancelled in $B_n/n$
\begin{equation} \label{eqn-adams-triv-div}
   \prod_{p-1 \notdiv n} p^{\ord_p n} \, \pdiv \, B_n \,.
\end{equation}
This is now known as Adams' theorem, although Adams \cite{adams78bern} only predicted
without proof that $p \pdiv n$ implies $p \pdiv B_n$ for primes $p-1 \notdiv n$ regarding
the table of the first 62 Bernoulli numbers $B_{2m}$ he had calculated in 1878.
On the other side, the property that $B_n/n$ is a $p$-integer for $p-1 \notdiv n$ is
necessarily needed to formulate the so-called Kummer congruences given by E.\,E.~Kummer
\cite{kummer50bern} earlier in 1851.
\medskip

These congruences and its generalizations are important properties of Bernoulli numbers
which lead to a $p$-adic view giving interesting information about $B_n/n$.
Let $\eulerphi$ be Euler's totient function. The Kummer congruences state
for $n, m, p, r \in \NN$, $n, m$ even, $p$ prime and $p-1 \notdiv n$
\begin{equation} \label{eqn-kummer-congr}
   (1-p^{n-1}) \frac{B_n}{n} \equiv (1-p^{m-1}) \frac{B_{m}}{m} \pmod{p^r}
\end{equation}
with $n \equiv m \ \mod{\eulerphi(p^r)}$.
\medskip

In 1850 Kummer \cite{kummer50flt} introduced the classification of regular and irregular
primes to characterize solutions of the famous Fermat's last theorem (FLT).
An odd prime $p$ is called \textsl{regular} if $p$ does not divide the class number of
the cyclotomic field $\QQ(\mu_p)$ with $\mu_p$ as the set of $p$-th roots of unity,
otherwise \textsl{irregular}. Kummer proved that if $p$ is regular then
FLT has no solution for the exponent $p$. He also
gave an equivalent definition concerning Bernoulli numbers:
An odd prime $p$ is called regular if $p$ does not divide any Bernoulli number
$B_n$ for $n=2,4,\ldots,p-3$, otherwise irregular. The \textsl{index of irregularity}
$i(p)$ counts these indices  for which $p \pdiv B_n$ happens.
In this case the pair $(p,n)$ is called an \textsl{irregular pair}.
First irregular primes are 37, 59, 67, 101.
\medskip

Regarding Bernoulli numbers, it will be very useful to combine properties of $B_n$
as well of $B_n/n$, the so-called divided Bernoulli number.
An easy consequence of the Kummer congruences provides that the numerator of
$B_n/n$ consists only of irregular primes and
that infinitely many irregular primes exist.
For the latter see a short proof of Carlitz \cite{carlitz54irrprime}, see also
\cite[Theorem 6, p.~241]{ireland90classical}.
Unfortunately, the more difficult question is still open
whether infinitely many regular primes exist. However,
calculations of Buhler, Crandall, Ernvall, Mets{\"a}nkyl{\"a}, and Shokrollahi
\cite{irrprime12M} show that about 60\% of all primes less
12 million are regular which agree with an expected distribution proposed by
Siegel \cite{siegel64irr}.
\medskip

All these basic results of Bernoulli numbers can be found in
the book of Ireland and Rosen \cite[Chapter 15]{ireland90classical}.
Throughout this paper all indices concerning Bernoulli numbers will be even and
$p$ an odd prime. Note that in older references the enumerating of Bernoulli numbers
can differ by a factor 2. Let $p^r \pdiveq n$ denote the highest power of $p$
dividing $n$ in order that $r=\ord_p n$.

%% file: sb_prelim.tex
\pagebreak

\section{Preliminaries}
\setcounter{equation}{0}

Here we will recall necessary facts about irregular prime powers
of Bernoulli numbers and $p$-adic zeta functions.
The definition of irregular pairs can be extended to irregular prime powers
which was already given in \cite{kellner04irrpowbn},
first introduced by the author \cite[Section 2.5]{kellner02irrpairord}.

\begin{defin} \label{def-ordn-irrpair}
A pair $(p,l)$ is called an \textsl{irregular pair of order $n$} if
$p^n \pdiv B_l/l$ with $2 \leq l < \eulerphi(p^n)$ and even $l$.
Let
\[
   \IRR_n := \{ (p, l) \sep p^n \pdiv B_l/l, \,\,
     2 \leq l < \eulerphi(p^n), \,\,  2 \pdiv l \}
\]
be the set of irregular pairs of order $n$. For a prime $p$
the \textsl{index} of irregular pairs of order $n$ is defined by
\[
    i_n(p) := \# \{ (p,l) \sep (p,l) \in \IRR_n \} \,.
\]
Let $(p,l) \in \IRR_n$ be an irregular pair of order $n$.
Let
\[
   (p,s_1,s_2, \ldots, s_n) \in \IRP_n \,, \quad
     l = \sum_{\nu=1}^{n} s_\nu \, \eulerphi( p^{\nu-1} )
\]
be the $p$-adic notation of $(p,l)$
with $0 \leq s_\nu < p$ for $\nu = 1,\ldots,n$ and $2 \pdiv s_1$, $2 \leq s_1 \leq p-3$.
The corresponding set will be denoted as $\IRP_n$. The pairs
$(p,l)$ and $(p,s_1,s_2, \ldots, s_n)$ will be called associated.
Let $(p,l) \in \IRR_n$ be an irregular pair of order $n$. Then define
\[
   \Delta_{(p,l)} \, \equiv \, p^{-n} \left(
     \frac{B_{l + \eulerphi(p^n)}}{l + \eulerphi(p^n)} - \frac{B_l}{l} \right) \pmod{p}
\]
with $0 \leq \Delta_{(p,l)} < p$.
In the case $\Delta_{(p,l)} = 0$ we will denote $\Delta_{(p,l)}$ as
\textsl{singular}.
\end{defin}
\medskip

\begin{remark} \label{rem-ordn-irrpair}
Note that this definition includes the usual definition of irregular pairs for $n=1$
with $i(p) = i_1(p)$. By Kummer congruences \refeqn{eqn-kummer-congr} the interval
$[2,\eulerphi(p^n)-2]$ is given for irregular pairs of order $n$ if they exist.
Moreover, we have the property that if $(p,l) \in \IRR_n$ then
$p^n \pdiv B_{l+k \eulerphi(p^n)}/(l+k \eulerphi(p^n))$ is valid
for all $k \in \NN_0$.
For simplification $(p,s_1,s_2, \ldots, s_n)$ is also called an irregular \textsl{pair}
with $(s_1,s_2, \ldots, s_n)$ as the second parameter in a $p$-adic manner.
It is easy to see that
if $(p,l_{n+1}) \in \IRR_{n+1}$ exists then $(p,l_n) \in \IRR_n$ exists with
$(p,l_n) = (p,l_{n+1} \ \mod{\eulerphi(p^n)})$, too.
\end{remark}
\medskip

The following proposition, see \cite[Prop.~5.3]{kellner04irrpowbn},
gives an unconditional representation of Bernoulli numbers by sets $\IRR_\nu$.
This is seen by \refeqn{eqn-clausen-von-staudt}, \refeqn{eqn-adams-triv-div},
and counting irregular prime powers.

\begin{prop} \label{prop-bn-repr}
Let $n$ be an even positive integer, then
\[
   B_n = (-1)^{\frac{n}{2}-1} \,
     \prod\limits_{p-1 \notdiv n} p^{\tau(p,n)+\ord_p n} \,\,\Big/\,\,
       \prod\limits_{p-1 \pdiv n} p
\]
with
\[
   \tau(p,n) := \sum_{\nu=1}^\infty \#( \, \IRR_\nu \cap
     \{ ( p, n \ \mod{\eulerphi(p^\nu)}) \} \, ) \,.
\]
\end{prop}
\medskip

The divided Bernoulli numbers $B_n/n$ are directly related to the Riemann zeta function
$\zeta(s)$ on the negative x-axis
\begin{equation} \label{eqn-zeta-1-n}
   \zeta( 1-n ) = - \frac{B_{n}}{n} \,, \qquad n \in \NN \,, \,\, n \geq 2 \,,
\end{equation}
where the Riemann zeta function is usually defined by the sum or the Euler product
\begin{equation} \label{eqn-zeta-sum-prod}
   \zeta(s) = \sum_{\nu=1}^\infty \nu^{-s}
      = \prod_{p} (1-p^{-s})^{-1} \,, \quad s \in \CC \,, \,\, \real s > 1 \,.
\end{equation}

Let $\ZZ_p$ be the ring of $p$-adic integers and $\QQ_p$ be the field of $p$-adic numbers.
Let $|\ |_p$ be the ultrametric absolute value on $\ZZ_p$ which is usually defined by
$|x|_p = p^{-\ord_p x}$. The property of Kummer congruences
allows the construction of a unique continuous $p$-adic zeta function which was introduced
by T.\,Kubota and H.\,W.\,Leopoldt \cite{leopoldt64zeta} in  1964,
see also Koblitz \cite[Chapter~II]{koblitz96padic}.
Modifying \refeqn{eqn-zeta-1-n} with an Euler factor $(1-p^{n-1})$ implies the
following definition.

\begin{defin} \label{def-zeta-padic}
Let $p$ be a prime with $p \geq 5$. Let
\[
   \zeta_p(1-n) \, := \, (1-p^{n-1}) \, \zeta(1-n)
     \, = \,  (1-p^{n-1}) \, \left( - \frac{B_{n}}{n} \right) \,.
\]
For a fixed $s_1 \in \{ 2,4,\ldots,p-3 \}$, define the
$p$-adic zeta function by
\[
   \zeta_{p,\,s_1}:\,\, \ZZ_p \to \ZZ_p \,, \qquad
     \zeta_{p,\,s_1} (s) := \lim_{t_\nu \to s} \,
       \zeta_p \big( 1-(s_1 +(p-1) t_\nu) \big)
\]
for $p$-adic integer $s$ by taking any sequence $(t_\nu)_{\nu \geq 1}$ of nonnegative integers
which $p$-adically converges to $s$.
\end{defin}
\medskip

By construction the $p$-adic zeta function $\zeta_{p,\,s_1} (s)$ interpolates
the zeta function $\zeta_p(1-n)$ at nonnegative integer values $s$ by
\[
    \zeta_{p,\,s_1}(s) = \zeta_p(1-n)
\]
with $n \equiv s_1 \ \mod{p-1}$ and $n=s_1+(p-1)s$.
Because of Kummer congruences
\[
   \zeta_p(1-n) \equiv \zeta_p(1-n') \pmod{p^r}
\]
for $n \equiv n' \ \mod{\eulerphi(p^r)}$ with $n \equiv n' \equiv s_1 \ \mod{p-1}$,
the $p$-adic zeta function $\zeta_{p,\,s_1} (s)$ is a unique continuous function on $\ZZ_p$
by means of interpolating property.

%% file: sb_delta.tex
\pagebreak

\section{\texorpdfstring{The $\Delta$-Conjecture}{The Delta-Conjecture}}
\setcounter{equation}{0}

Let $(p,l) \in \IRR_1$ be an irregular pair, then we are interested in the
behavior of $\Delta_{(p,l)}$. Essentially, there are two cases to consider:
$\Delta_{(p,l)} \neq 0$ and $\Delta_{(p,l)} = 0$.
Now, calculations in \cite{irrprime12M} for irregular primes $p < 12\,000\,000$ show
that $\Delta_{(p,l)} \neq 0$ is always valid. No singular $\Delta_{(p,l)}$ has been
found yet. However, the improbable case of a singular $\Delta_{(p,l)}$
which implies a strange behavior without regularity is described in the next section.
The following theorem gives the main result of irregular pairs in the nonsingular case,
see \cite[Theorem~3.1]{kellner04irrpowbn}.

\begin{theorem} \label{theor-delta-irr}
Let $(p,l_1) \in \IRR_1$ be an irregular pair with $\Delta_{(p,l_1)} \neq 0$.
Then for each \mbox{$n > 1$} there exists exactly one irregular pair of order $n$
corresponding to $(p,l_1)$. Therefore, a unique sequence
$(l_n)_{n \geq 1}$ resp. $(s_n)_{n \geq 1}$ exists with
\[
   (p,l_n) \in \IRR_n \quad \mbox{resp.} \quad (p,s_1,\ldots,s_n) \in \IRP_n
\]
and
\[
    l_1 \leq l_2 \leq l_3 \leq \ldots \,,
      \quad \lim_{n \to \infty} l_n = \infty \,.
\]
Moreover, one has
\[
   \Delta_{(p,l_1)} = \Delta_{(p,l_2)} = \Delta_{(p,l_3)} = \ldots \,\,.
\]
If $\Delta_{(p,l_{1,\nu)}} \neq 0$ for all $i(p)$ irregular pairs $(p,l_{1,\nu}) \in \IRR_1$, then
\[
   i(p) = i_2(p) = i_3(p) = \ldots \,\,.
\]
\end{theorem}

\begin{defin} \label{def-char-padic}
Let $(p,l) \in \IRR_1$ be an irregular pair with $\Delta_{(p,l)} \neq 0$.
Then Theorem \ref{theor-delta-irr} provides a unique sequence
$(s_\nu)_{\nu \geq 1}$ with $l=s_1$.
Define a \textsl{characteristic} $p$-adic integer
\[
   \chi_{(p,\,l)} = \sum_{\nu \geq 0} s_{\nu+2} \, p^\nu \,\, \in \ZZ_p
\]
which contains all information of irregular pairs of higher order
corresponding to $(p,l)$.
\end{defin}
\medskip

The following theorem, a result of \cite[Theorem~4.6/4.10]{kellner04irrpowbn},
shows the behavior of the $p$-adic zeta function.
The Kummer congruences \refeqn{eqn-kummer-congr} are valid
by the implication
\[
   n \equiv m \pmod{\eulerphi(p^r)} \quad \Longrightarrow \quad
     (1-p^{n-1}) \frac{B_n}{n} \equiv (1-p^{m-1}) \frac{B_{m}}{m} \pmod{p^r} \,,
\]
but the converse does not hold in general. The first nontrivial counterexample is given by $p=13$ and
$B_{16}/16 - B_4/4 = -7\cdot13^2/2720$. Note that also $B_{14}/14 - B_2/2 = 0$ happens which
is the only exception that divided Bernoulli numbers are equal.

\begin{theorem} \label{theor-zeta-padic-zero}
Let $(p,l) \in \IRR_1$ be an irregular pair with $\Delta_{(p,l)} \neq 0$.
Let $s, t \in \ZZ_p$.
Then the $p$-adic zeta function $\zeta_{p,\,l}$ has a unique zero with
$\zeta_{p,\,l} \!\left( \chi_{(p,\,l)} \right) = 0$.
A strong version of the Kummer congruences holds
\[
   |\zeta_{p,\,l} (s) - \zeta_{p,\,l} (t)|_p = |p \, ( s - t ) |_p \,.
\]
Moreover, one has
\[
   \frac{\zeta_{p,\,l} (s) - \zeta_{p,\,l} (t)}{p \, ( s - t )}
      \equiv - \Delta_{(p,l)} \pmod{p\ZZ_p}
         \,, \quad s \neq t \,.
\]
\end{theorem}

As a consequence, we easily obtain
$|\zeta_{p,\,l} (s)|_p = |p \, ( \chi_{(p,l)} - s ) |_p$ for $s \in \ZZ_p$
under the assumption above. Thus, the power of the irregular prime $p$ can be
described by measuring the $p$-adic distance to the zero $\chi_{(p,\,l)}$ of
the $p$-adic zeta function $\zeta_{p,\,l}$.
Since $|p \, n|_p \, |\zeta(1-n)|_p = 1$ for $n \equiv 0 \ \mod{p-1}$ and
$|\zeta(1-n)|_p = |\zeta_{p,\,l} (s)|_p$ for $s = \frac{n-l}{p-1}$, $n \equiv l$ $\mod{p-1}$,
we obtain a structural formula of $\zeta(1-n)$ resp.\ $B_n/n$, see \cite[Theorem~4.9]{kellner04irrpowbn}.
Note that we combine the numerator and denominator of $B_n/n$ in this formula.
Vaguely speaking, the numerator can be described by zeros of $\zeta_{p,\,l}$ and
the denominator by poles of $\zeta_{p,\,0}$ lying at 0, where $\zeta_{p,\,0}:\, \ZZ_p \to \QQ_p$
extends Definition \ref{def-zeta-padic} using arguments given in \cite[Chapter II, p.~46]{koblitz96padic}.

\begin{theorem} \label{theor-zeta-prod-psi-0}
Let $\PP$ be the set of primes. Then define $\Psi_0 = \IRR_1 \cup ( \PP \times \{0\} )$
and $\chi_{(p,0)}=0$ for all $p \in \PP$. Define $\rho(l)=1-2\sign(l)=\pm1$ for $l \geq 0$.
Let $n$ be an even positive integer, then
under the assumption that no singular $\Delta_{(p,l)}$ exists, one has
\[
   \zeta(1-n) = (-1)^{\frac{n}{2}} \!\!\!
     \prod\limits_{(p,l)\in\Psi_0 \atop n \equiv l \,\mods{p-1}}
       \left( \frac{\,|\chi_{(p,l)} - \frac{n-l}{p-1}|_p}{p} \right)^{\rho(l)} \,.
\]
\end{theorem}

One may ask whether the structure of the Riemann zeta function at odd negative integers is given
by this simple form. Now, all these facts substantiated by calculations lead
to the following conjecture, already mentioned in \cite{kellner02irrpairord}.

\begin{conj}[$\Delta$-Conjecture] \label{conj-delta}
Let $p$ be an irregular prime. Then the following properties hold:
\begin{enumerate}
\item $\Delta_{(p,l_\nu)}$ is not singular for all irregular pairs $(p,l_\nu) \in \IRR_1$.
\item $i(p)=i_2(p)=i_3(p)=\ldots\,$.
\item The $p$-adic zeta function $\zeta_{p,\,l_\nu}$ has a unique zero $\chi_{(p,\,l_\nu)}$.
\item A strong form of the Kummer congruences holds for $\zeta_{p,\,l_\nu}$.
\end{enumerate}
\end{conj}

Finally, assuming the $\Delta$-Con\-jecture, we also obtain a structural formula of
Bernoulli numbers which gives a decomposition in three products. The first and last
product are trivially given, the complicated product in the middle consists only of
irregular primes.

\begin{theorem} \label{theor-bn-prod-3}
Let $n$ be an even positive integer, then under the assumption of the $\Delta$-Conjecture
\[
   B_n = (-1)^{\frac{n}{2}-1} \prod_{p-1 \notdiv n} |n|_p^{-1}
     \!\! \prod\limits_{(p,l)\in\IRR_1 \atop n \equiv l \,\mods{p-1}}
       \!\!\!\!\! |p \, (\chi_{(p,l)} - {\textstyle \frac{n-l}{p-1}})|_p^{-1}
          \,\,\, \prod\limits_{p-1 \pdiv n} p^{-1} \,.
\]
\end{theorem}

\begin{proof}
The product in the middle is given by Theorem \ref{theor-zeta-prod-psi-0}. The
first resp.\ last product is a consequence of \refeqn{eqn-adams-triv-div} resp.\
\refeqn{eqn-clausen-von-staudt}.
\end{proof}

The main task remains to determine the zero of a $p$-adic zeta function associated
with an irregular pair $(p,l)$. An irregular pair $(p,l_n) \in \IRR_n$ of order $n$
yields an approximation of the zero $\chi_{(p,\,l)}$. Fortunately, these
irregular pairs of higher order can be computed with little effort by calculating
a small number of divided Bernoulli numbers with relatively small indices.
For algorithms and calculated pairs see \cite[Section~5,Table~A.3]{kellner04irrpowbn}.
For example, we have
\[
    (157,62,40,145,67,29,69,0,87,89,21) \in \IRP_{10}
\]
which also shows a rare occurrence of a zero in the $p$-adic sequence, here $s_7=0$.
This is the only zero which occurs in the $p$-adic sequence of irregular pairs
of order 10 for $p < 1000$, see \cite[Table~A.3]{kellner04irrpowbn}. This means
\[
    (157, 6\,557\,686\,520\,486) \in \IRR_6 \cap \IRR_7
\]
respectively, since the index is prime to 157 and $s_8=87$,
\[
    157^7 \pdiveq B_{6\,557\,686\,520\,486} \,.
\]
Now, we will never be able to compute this giant Bernoulli number nor
we can use Kummer congruences to verify this result directly!

%% file: sb_sing.tex
\section{The singular case}
\setcounter{equation}{0}

The case of a singular $\Delta_{(p,l)}$ is more complicated than the nonsingular
case where we have a certain regularity.
However, no such singular $\Delta_{(p,l)}$ has been found yet. Therefore,
supported by massive computations, we can regard such an event as a rare exception.
The following theorem, see \cite[Theorem~3.2]{kellner04irrpowbn},
shows the strange behavior of irregular pairs of higher order with a
singular $\Delta_{(p,l)}$.

\begin{theorem} \label{theor-delta-irr-0}
Let $(p,l_n) \in \IRR_n$ be an irregular pair of order $n$
with $\Delta_{(p,l_n)}=0$. Then there exist two cases:
\begin{enumerate}
 \item $(p,l_n) \notin \IRR_{n+1}:$
   There are no irregular pairs of order $n+1$ and higher.
 \item $(p,l_n) \in \IRR_{n+1}:$
   There exist $p$ irregular pairs
   $(p,l_{n+1,j}) = (p,l_n+j \eulerphi(p^n)) \in \IRR_{n+1}$ of order $n+1$
   with $\Delta_{(p, \,l_{n+1,j})} = 0$ for $j=0,\ldots,p-1$.
\end{enumerate}
\end{theorem}

The following diagram demonstrates the behavior of the singular case.
This situation can be described by a rooted $p$-ary tree
of irregular pairs of higher order.

\unitlength0.0075cm
\begin{center}
\begin{picture}(700,520)
  \put(200,460){$\Delta_{(p,l_1)} = 0$}
  \put(10,390){$\IRR_1$}
  \put(10,290){$\IRR_2$}
  \put(10,190){$\IRR_3$}
  \put(10,090){$\IRR_4$}
  \put(145,030){$l_1$}
  \put(245,030){$l_2$}
  \put(395,030){$l_3$}
  \put(150,400){\line(0,-1){100}}
  \put(150,400){\line(1,-1){100}}
  \put(150,400){\line(2,-1){200}}
  \put(150,400){\circle*{10}}
  \put(150,300){\circle*{10}}
  \put(250,300){\circle*{10}}
  \put(350,300){\circle*{10}}

  \put(250,300){\line(0,-1){100}}
  \put(250,300){\line(3,-4){75}}
  \put(250,300){\line(3,-2){150}}

  \put(325,200){\circle*{10}}
  \put(250,200){\circle*{10}}
  \put(400,200){\circle*{10}}

  \put(250,200){\line(0,-1){100}}
  \put(250,200){\line(1,-2){50}}
  \put(250,200){\line(1,-1){100}}
  \put(400,200){\line(0,-1){100}}
  \put(400,200){\line(1,-2){50}}
  \put(400,200){\line(1,-1){100}}

  \put(250,100){\circle*{10}}
  \put(300,100){\circle*{10}}
  \put(350,100){\circle*{10}}
  \put(400,100){\circle*{10}}
  \put(450,100){\circle*{10}}
  \put(500,100){\circle*{10}}
\end{picture}
\end{center}

Here a vertical line indicates that $(p,l_n) \in \IRR_n \cap \IRR_{n+1}$ happens.
We then have $p$ irregular pairs of order $n+1$ which are
represented by branches. In this case, the corresponding Bernoulli number $B_{l_n}/l_n$
\textsl{decides} whether there exist further branches or they stop.
Instead of $n$ the order of the $p$-power must be at least $n+1$.
This also means that an associated irregular pair $(p,s_1,\ldots,s_{n+1}) \in \IRP_{n+1}$
must have a zero $s_{n+1}=0$ in its $p$-adic notation each time.
Now, it is worth saying that no irregular pair
$(p,l)$ has been found with $p^2 \pdiv B_l$ resp.\ $(p,l) \in \IRR_1 \cap \IRR_2$
for $p < 12\,000\,000$, see \cite{irrprime12M}, while an example of an element of
$\IRR_6 \cap \IRR_7$ is shown in the previous section.

\begin{defin}
Let $(p,l) \in \IRR_1$ be an irregular pair with a singular $\Delta_{(p,l)}$.
Then define a rooted $p$-ary tree of irregular pairs of higher order like in the
diagram above given by Theorem \ref{theor-delta-irr-0}.
Each node contains one irregular pair of higher order. Note that these pairs are not
necessarily distinct. We denote this tree as $T^0_{(p,l)}$ related to the root node $(p,l)$.
The tree $T^0_{(p,l)}$ has the property that each node of height $r$ lies in $\IRR_{r+1}$.
A tree $T^0_{(p,l)} = \{ (p,l) \}$ is called a \textsl{trivial tree} having height 0.
A tree of height one is given by the root node $(p,l)$ and its $p$ child nodes
$(p,l+j\eulerphi(p))$ with $j=0,\ldots,p-1$. A tree with height $\geq 2$ always contains the
latter one.
\end{defin}
\medskip

In the nonsingular case, we have a zero of the $p$-adic zeta function. In contrast to,
the singular case does not guarantee that irregular pairs of higher order exist at all.
However, an exception does not destroy Theorem \ref{theor-bn-prod-3}
but complicates the formula, because we then have to consider the tree $T^0_{(p,l)}$
of irregular pairs of higher order. Thus, we obtain an
unconditional formula by combining both cases.

\begin{theorem} \label{theor-bn-prod-4}
Let $n$ be an even positive integer, then
\begin{eqnarray*}
  B_n &=& (-1)^{\frac{n}{2}-1} \prod_{p-1 \notdiv n} |n|_p^{-1}
    \!\! \prod\limits_{(p,l)\in\IRR_1, \,\Delta_{(p,l)} \neq 0 \atop n \equiv l \,\mods{p-1}}
      \!\!\!\!\! |p \, (\chi_{(p,l)} - {\textstyle \frac{n-l}{p-1}})|_p^{-1} \\
  && \quad \prod\limits_{(p,l)\in\IRR_1, \,\Delta_{(p,l)} = 0 \atop n \equiv l \,\mods{p-1}}
      \!\!\!\!\! p^{1+h_0(p,n)}
        \,\,\, \prod\limits_{p-1 \pdiv n} p^{-1}
\end{eqnarray*}
with the height of $(p,n)$ defined by
\[
   h_0(p,n) = \max \left\{ \mathrm{height}( (p,l') ) \sep
     (p,l') \in T^0_{(p,l)} \cap \{ (p, n \ \mod{\eulerphi(p^\nu)}) \}_{\nu \geq 1}
       \right\} \,.
\]
Moreover, $h_0(p,n) = 0$ $\iff$ the tree $T^0_{(p,l)}$ is trivial. Otherwise $h_0(p,n) \geq 1$.
\end{theorem}

\begin{proof}
The case $\Delta_{(p,l)} \neq 0$ is already handled by Theorem \ref{theor-bn-prod-3}.
Now, assume $\Delta_{(p,l)} = 0$ with a given tree $T^0_{(p,l)}$.
We have to determine the max.\ height of a node $(p,l_{\nu,j}) \in T^0_{(p,l)} \cap \IRR_\nu$
which is equal to $( p, n \ \mod{\eulerphi(p^\nu)})$, as a consequence of Remark \ref{rem-ordn-irrpair}.
The root node $(p,l)$ has height 0, so the exponent equals $1 + h_0(p,n)$.
If the tree $T^0_{(p,l)}$ is trivial, then $h_0(p,n) = 0$ is constant.
On the other side, a tree $T^0_{(p,l)}$ having height $\geq 1$ contains $p$ irregular pairs
of order two. Then $(p, n \ \mod{\eulerphi(p^2)}) \in T^0_{(p,l)}$ is always valid
which finally yields $h_0(p,n) \geq 1$.
\end{proof}

\begin{corl}
Let $n$ be an even positive integer, then
\[
  \zeta(1-n) = (-1)^{\frac{n}{2}}
     \!\! \prod\limits_{(p,l)\in\IRR_1, \,\Delta_{(p,l)} \neq 0 \atop n \equiv l \,\mods{p-1}}
     \!\!\!\!\! |p \, (\chi_{(p,l)} - {\textstyle \frac{n-l}{p-1}})|_p^{-1}
     \!\!\! \prod\limits_{(p,l)\in\IRR_1, \,\Delta_{(p,l)} = 0 \atop n \equiv l \,\mods{p-1}}
      \!\!\!\!\! p^{1+h_0(p,n)}
        \,\,\, \prod\limits_{p-1 \pdiv n} \frac{|n|_p}{p}
\]
with $h_0(p,n)$ as defined above.
\end{corl}

%% file: sb_appl.tex
\section{Applications}
\setcounter{equation}{0}

Regarding Theorem \ref{theor-bn-prod-4} and the definitions of $h_0$ and $\chi_{(p,l)}$,
we can state an extended version of Adams' theorem given by \refeqn{eqn-adams-triv-div}.

\begin{theorem} \label{theor-ext-adam}
Let $n$ be an even positive integer. Let $p$ be a prime with $p^r \pdiveq n$,
$r \geq 1$, and $p-1 \notdiv n$. Let $l \equiv n \ \mod{p-1}$ with $0 < l < p-1$.
Then $p^{r+\delta} \pdiveq B_n$ with the following cases:
\begin{enumerate}
\item If $p$ is regular, then $\delta = 0$.
\item If $p$ is irregular with $(p,l)\notin\IRR_1$, then $\delta = 0$.
\item If $p$ is irregular with $(p,l)\in\IRR_1$, $\Delta_{(p,l)} \neq 0$,
      then $\delta = 1 + \ord_p \, (\chi_{(p,l)} - {\textstyle \frac{n-l}{p-1}})$.
\item If $p$ is irregular with $(p,l)\in\IRR_1$, $\Delta_{(p,l)} = 0$,
      then $\delta = 1 + h_0(p,n)$.
\end{enumerate}
Additionally, in case (3) resp.\ (4), if $(p,l,l) \notin \IRP_2$,
then $\delta = 1$ is bounded, otherwise $\delta \geq 2$.
\end{theorem}

\begin{proof}
We have to consider the formula of Theorem \ref{theor-bn-prod-4},
then the first product yields $p^r \pdiv B_n$. Only the second resp.\ third product
can give additional $p$-factors. Therefore, case (1) and (2) are given by definition.
Now, we can assume $(p,l)\in\IRR_1$. A nonsingular $\Delta_{(p,l)}$ provides
$\delta = \ord_p |p \, (\chi_{(p,l)} - {\textstyle \frac{n-l}{p-1}})|_p^{-1} =
   1 + \ord_p \, (\chi_{(p,l)} - {\textstyle \frac{n-l}{p-1}})$ in case (3).
On the other side, a singular $\Delta_{(p,l)}$ provides $\delta = 1 + h_0(p,n)$ in case (4).
The additional cases are shown as follows.
\smallskip

Case (3): By assumption, $n=p^r n'$ with some integer $n'$. We have to evaluate
\[
   d = \ord_p \, (\chi_{(p,l)} - {\textstyle \frac{n-l}{p-1}}) =
     \ord_p \, (p \,\chi_{(p,l)} - \chi_{(p,l)} + l - p^r n') \,.
\]
Since $r \geq 1$, we $p$-adically obtain $\chi_{(p,l)} = l + s_3\,p + \ldots \iff d \geq 1$
which is equal to $(p,l,l) \in \IRP_2$. Conversely, $(p,l,l) \notin \IRP_2$ yields $d=0$.
Case (4): $(p,l,l) \notin \IRP_2$ $\iff$ the tree $T^0_{(p,l)}$ is trivial.
Then $h_0(p,n)=0$ is constant. Conversely, $(p,l,l) \in \IRP_2$ yields $h_0(p,n) \geq 1$
as a result of Theorem \ref{theor-bn-prod-4}.
\end{proof}
\medskip

So far, no $(p,l,l) \in \IRP_2$ has been found yet. We can even raise the
value $\delta$ in the following way.

\begin{corl}
Assume that $(p,l,\ldots,l) \in \IRP_r$ exists with some $r\geq1$. Let $n=lp^r$. Then, we have
$p^r \pdiveq n$ and $p^{2r} \pdiv B_n$.
\end{corl}

\begin{proof}
By Definition \ref{def-ordn-irrpair}, we have
\[
   l_r = \sum_{\nu=1}^{r} l \, \eulerphi( p^{\nu-1} ) = l p^{r-1}
     \quad \mbox{with} \quad (p,l_r) \in \IRR_r \,.
\]
Then $p^r \pdiv B_{l_r+k \, \eulerphi(p^r)}/(l_r+k \, \eulerphi(p^r))$ is valid for all $k\geq0$.
Choose $n=l_r + l\,\eulerphi(p^r)=lp^r$. Thus, $p^r \pdiv B_n/n$ and finally $p^{2r} \pdiv B_n$.
Note that we cannot predict that $p^{2r} \pdiveq B_n$ in general.
\end{proof}

Johnson \cite[Theorem, p.~655]{Johnson74} calculated irregular pairs up to $p < 8000$.
Correspondingly, he also calculated the now called irregular pairs $(p,s_1,s_2) \in \IRP_2$
of order two in that range, proving that $(p,l,l) \notin \IRP_2$ for $p < 8000$. See also
\cite[Table~A.3]{kellner04irrpowbn} for calculations of irregular pairs of order 10 for $p < 1000$.
In a similar manner, the nonexistence of irregular pairs $(p,l,l-1)$ of order two
plays an important role in Iwasawa theory, see Washington \cite{washington97cyclo} for
Iwasawa theory and \cite[Section~6]{kellner04irrpowbn} for this special result.
In context of cyclotomic invariants,
calculations of \cite{irrprime12M} ensure that no $(p,l,l-1) \in \IRP_2$ exists for
$p < 12\,000\,000$. One may conjecture that no such special irregular pairs
$(p,l,l)$ resp.\ $(p,l,l-1)$ of order two exist. But there is still a long way
to prove such results, even to understand properly which role the zeros $\chi_{(p,l)}$ play.
\medskip

By Definition \ref{def-ordn-irrpair}, we have the relation
\[
   (p,l,l) \notin \IRP_2  \iff  p^2 \notdiv B_{lp}/(lp)
     \iff  p^3 \notdiv B_{lp} \,.
\]
Yamaguchi \cite{Yama76bern} also verified by calculation that $p^3 \notdiv B_{lp}$
for all irregular pairs $(p,l)$ with $p < 5500$, noting
that this was conjectured earlier by Morishima in general.
Interestingly, the condition $p^3 \notdiv B_{lp}$ is related to the second
case of FLT, see \cite[Theorem~9.4, p.~174]{washington97cyclo}. See also
\cite[Corollary~8.23, p.~162]{washington97cyclo} for a different context.
Under the assumption of the conjecture of Kummer--Vandiver and that
no $(p,l,l) \in \IRP_2$ exists, the second case of FLT has no solution for the
exponent $p$. For details we refer to the references cited above.

\begin{remark}
Thangadurai \cite{thang04adam} also investigates Adams' theorem
claiming a conjecture that $\delta \in \{0,1\}$ where $\delta$ is given
as in Theorem \ref{theor-ext-adam}. The following theorem \cite[Theorem~2.5, p.~172]{thang04adam}
is formulated, here translated into our terminology
\[
   \delta \in \{0,1\} \iff p^3 \notdiv B_{lp} \iff p^2 \notdiv B_l \quad \mbox{for all \ }
     l = 2,4,\ldots,p-3 \,.
\]
The first equivalence agrees with our results, but the second equivalence is \textbf{false}.
Theorem 2.3 and Theorem 2.4 of \cite[pp.~171--172]{thang04adam} are incorrect in general.
Theorem 2.3 of \cite{thang04adam} is an incorrect citation of the results of
Johnson \cite[Theorem, p.~655]{Johnson74} which are only valid for
irregular primes $p < 8000$ and were verified by calculations.
In our words, Theorem 2.3 of \cite{thang04adam} would imply
$\Delta_{(p,l)} \neq 0 \Longrightarrow (p,l,l) \notin \IRP_2$ for all irregular pairs!
But $\Delta_{(p,l)} \neq 0$ only implies the existence of a $(p,l,k) \in \IRP_2$
with a unique $k$ in the range $0 \leq k < p$. Theorem 2.4 of \cite{thang04adam} states
the following: Let $(p,l) \in \IRR_1$. If $p \pdiveq B_{l+k(p-1)}/(l+k(p-1))$
for some $k \geq 1$, then $p \pdiveq B_l$.
Now, this theorem would imply that $p^2 \notdiv B_l$ in the nonsingular case! Counterexample:
$\Delta_{(p,l)} \neq 0 \Longrightarrow$ there exists exactly one $(p,l,k') \in \IRP_2$.
The case $k'=0$ can occur. Then we have $p^2 \pdiv B_l$, but
$p \pdiveq B_{l+(p-1)}/(l+(p-1))$. As a consequence,
\cite[Theorem~1, p.~170]{thang04adam} is only valid for irregular primes $p<8000$
by results of Johnson. Finally, \cite[Theorem~4.2, p.~176]{thang04adam}
is also incorrect, because $(p,l,l-1) \notin \IRP_2$ does not imply $p^2 \notdiv B_l$ in the
nonsingular case $\Delta_{(p,l)} \neq 0$, while in the singular case this implication always
holds as seen by a trivial tree $T^0_{(p,l)}$.
In conclusion, all depending results of these theorems cited above are incorrect, especially
the claimed equivalences to other hypotheses of Bernoulli numbers in \cite{thang04adam}.
\end{remark}
\medskip

Certainly, the converse of Adams' theorem does not hold, but one can state a somewhat different
result which deals with the common prime factors of numerators and denominators of
adjoining Bernoulli numbers, see \cite[Satz 2.3.4, p.~35]{kellner02irrpairord}.

\begin{theorem} \label{theor-bn-num-denom}
Let $\mathcal{S} = \{ 2,4,6,8,10,14 \}$ be the set of all even indices $m$ where the numerator
of $B_m/m$ equals 1. Write $B_n = \Bnum_n/\Bden_n$ with $(\Bnum_n,\Bden_n)=1$.
Let $k, n$ be even positive integers with $k \in \mathcal{S}$ and $n-k \geq 2$. Then
\[
    D=(\Bnum_n,\Bden_{n-k}) \quad \mbox{provides} \quad D \pdiv n \,.
\]
Moreover, if $D > 1$ then $D=p_1 \cdots p_r$ with $r \geq 1$ and $p_\nu \notdiv \Bden_k$,
$p_\nu \notdiv B_n/n$ for all $\nu=1,\ldots,r$.
\end{theorem}

\begin{proof}
Assume $D > 1$. We then have $D=p_1 \cdots p_r$ with $r \geq 1$ since $\Bden_{n-k}$ is squarefree
by \refeqn{eqn-clausen-von-staudt}.
Let $\nu \in \{ 1, \ldots, r \}$. We have the following properties: $p_\nu \pdiv \Bnum_n$ and
$p_\nu - 1 \notdiv n$, additionally, $p_\nu \pdiv \Bden_{n-k}$ and $p_\nu - 1 \pdiv n-k$ with
$p_\nu < n$. Hence, we obtain $p_\nu \notdiv \Bden_k$. Assume to the contrary that
$p_\nu-1 \pdiv k$. By $p_\nu-1 \pdiv n-k$ we get $p_\nu-1 \pdiv k + n-k = n$. Contradiction.
\smallskip

Assume $p_\nu \notdiv n$ or $p_\nu \pdiv B_n/n$,
then we obtain by Kummer congruences \refeqn{eqn-kummer-congr}
\[
   0 \equiv \frac{B_n}{n} \equiv \frac{B_k}{k} \pmod{p_\nu} \,,
\]
since $n-k \equiv 0 \ \mod{p_\nu-1}$. By properties of the set $\mathcal{S}$
\begin{equation} \label{eqn-loc-bn-num-denom-1}
   \frac{B_k}{k} \not\equiv 0 \pmod{p_\nu}
\end{equation}
yields a contradiction. We obtain $p_\nu \pdiv n$ and $p_\nu \notdiv B_n/n$.
Finally $D \pdiv n$ is valid.
\end{proof}

Now, the set $\mathcal{S}$ cannot be enlarged, because \refeqn{eqn-loc-bn-num-denom-1}
does not hold in general for numerators having prime factors. For example, let $p=691$ and
$n=12+(p-1)$, then we have $p \pdiv B_{12}/12$ and
$D=(\Bnum_n,\Bden_{n-12})=pc \notdiv n$ with some $c \geq 1$.
Actually, $c=1$ with the help of \textbf{Mathematica}.
On the other hand, one trivially obtains for $k \in \mathcal{S}$,
$p$ prime with $p-1 \notdiv k$, $n=kp$ infinitely many examples of $D > 1$. In the following
theorem, Theorem \ref{theor-bn-num-denom} plays a crucial role.
Define for positive integers $n$ and $m$ the summation formula of consecutive integer powers by
\[
   S_n(m) = \sum_{\nu=0}^{m-1} \nu^n \,.
\]
Many congruences concerning function $S_n$ are naturally related to Bernoulli numbers.

\begin{theorem} \label{theor-bn-sn}
Let $n, m$ be positive integers with even $n$. Then
\[
   m^{r+1} \pdiv S_n(m) \quad \iff  \quad m^{r} \pdiv B_n
\]
for $r=1,\,2$.
\end{theorem}

\begin{proof}
Write $B_n = \Bnum_n/\Bden_n$ with $(\Bnum_n,\Bden_n)=1$.
Assume $m > 1$ and $n \geq 10$ with even $n$, otherwise we have $|\Bnum_n|=1$ for $n=2,4,6,8$.
It is well known, see \cite[p.~234]{ireland90classical}, that for even $n \geq 10$
\begin{equation} \label{eqn-loc-bn-sn-1}
  S_n(m) = B_n \, m + \binom{n}{2} B_{n-2} \frac{m^3}{3} +
    \sum_{k=4}^n \binom{n}{k} B_{n-k} \frac{m^{k+1}}{k+1} \,.
\end{equation}
We have to examine carefully the sum given in \refeqn{eqn-loc-bn-sn-1}.
By the result of Clausen--von Staudt \refeqn{eqn-clausen-von-staudt} the denominator
of all nonzero Bernoulli numbers is squarefree including $B_0$ and $B_1$.
Consider each prime factor $p^s \pdiveq m$. Then, we have
\begin{equation} \label{eqn-loc-bn-sn-2}
   \ord_p \left( \binom{n}{k} B_{n-k} \frac{m^{k+1}}{k+1} \right)
     \geq s(k+1) - 1 - \ord_p (k+1) \geq \left\{
     \begin{array}{rl}
        s \,, & \!\! k \geq 2, \, p \geq 2 \\
       2s \,, & \!\! k \geq 2, \, p \geq 5 \\
       3s \,, & \!\! k \geq 4, \, p \geq 5
     \end{array} \right.
\end{equation}
for such $k \leq n$ where $B_{n-k} \neq 0$. Critical cases are to consider for $p=2,3,5$ and
$s=1$. Now, we are ready to evaluate \refeqn{eqn-loc-bn-sn-1} $\mod{m^t}$ for certain $t$.
\medskip

\textbf{Case} $r=1$: Assume $(m,\Bden_n) > 1$, then we obtain by \refeqn{eqn-loc-bn-sn-2}
(case $k \geq 2, \, p \geq 2$)
\[
   S_n(m) \equiv B_n \, m \equiv \frac{\Bnum_n}{\Bden_n} \, m \not\equiv 0 \pmod{m} \,.
\]
Therefore, $(m,\Bden_n) = 1$, $2 \notdiv m$, $3 \notdiv m$, and $p \geq 5$ must hold. Hence,
by \refeqn{eqn-loc-bn-sn-2} (case $k \geq 2, \, p \geq 5$), we can write
$S_n(m) \equiv B_n \, m \ \mod{m^2}$. Consequently,
\[
   0 \equiv S_n(m) \equiv B_n \, m \pmod{m^2}
\]
provides $m^2 \pdiv S_n(m) \iff m \pdiv B_n$.
\medskip

\textbf{Case} $r=2$: We have $m \pdiv B_n$ and $(m,6)=1$, because either $m^2 \pdiv B_n$
or $m^3 \pdiv S_n(m)$ is assumed. The latter case implies $m^2 \pdiv S_n(m)$ and therefore
with case $r=1$ also $m \pdiv B_n$. By \refeqn{eqn-loc-bn-sn-2} (case $k \geq 4, \, p \geq 5$)
we obtain
\begin{equation} \label{eqn-loc-bn-sn-3}
   S_n(m) \equiv B_n \, m + \binom{n}{2} \frac{B_{n-2}}{3} \, m^3
     \equiv B_n \, m + \frac{n(n-1) \Bnum_{n-2}}{6 \Bden_{n-2}} \, m^3 \pmod{m^3} \,.
\end{equation}
Our goal is to show that the second term vanishes, but the denominator $\Bden_{n-2}$
could possibly remove prime factors from $m$. Now, Theorem \ref{theor-bn-num-denom} asserts that
$(\Bnum_n, \Bden_{n-2}) \pdiv n$. We also have $(m, \Bden_{n-2}) \pdiv n$ since $m \pdiv B_n$.
This means that the factor $n$ adds those primes which $\Bden_{n-2}$ possibly removes from $m$.
Therefore, the second term of \refeqn{eqn-loc-bn-sn-3} vanishes $\mod{m^3}$. The rest follows
again by $S_n(m) \equiv B_n \, m \equiv 0 \ \mod{m^3}$.
\end{proof}

One can improve the value $r$ for certain $m \pdiv n$, since $\binom{n}{k}$ appears in the sum
\refeqn{eqn-loc-bn-sn-1}, but not in general. Let $p=37$ and $l=37580$.
We then have $(p,l) \in \IRR_3$ and $p^3 \pdiv B_l$, but $p^4 \notdiv S_l(p)$
which was checked with \textbf{Mathematica}.
\pagebreak

\textbf{Example}:
\begin{enumerate}
\item We have $B_{42} = 1520097643918070802691/1806$.
Since the numerator $\Bnum_{42}$ is a big prime, we obtain for $m>1$
\[
   m^2 \pdiv S_{42}(m) \quad \Longleftrightarrow \quad m = 1520097643918070802691 \,.
\]
\item We have $\Bnum_{50} = 5^2 \cdot 417202699 \cdot 47464429777438199$ and
$\Bden_{48} = 2 \cdot 3 \cdot 5 \cdot 7 \cdot 13 \cdot 17$. Hence, for $m>1$
\[
   m^3 \pdiv S_{50}(m) \quad \Longleftrightarrow \quad m = 5 \,.
\]
\end{enumerate}

%% file: sb_main.bbl
\newcommand{\etalchar}[1]{$^{#1}$}
\begin{thebibliography}{BCE{\etalchar{+}}01}

\bibitem[Ada78]{adams78bern}
J.~C. Adams.
\newblock Table of the values of the first sixty-two numbers of {B}ernoulli.
\newblock {\em J. Reine Angew. Math.}, 85:269--272, 1878.

\bibitem[BCE{\etalchar{+}}01]{irrprime12M}
J.~Buhler, R.~Crandall, R.~Ernvall, T.~Mets{\"a}nkyl{\"a}, and M.~A.
  Shokrollahi.
\newblock Irregular primes and cyclotomic invariants to 12 million.
\newblock {\em Journal of Symbolic Computation}, 31(1/2):89--96, January 2001.

\bibitem[Car54]{carlitz54irrprime}
L.~Carlitz.
\newblock A note on irregular primes.
\newblock {\em Proc. Amer. Math. Soc.}, 5:329--331, 1954.

\bibitem[Cla40]{clausen40bern}
T.~Clausen.
\newblock Lehrsatz aus einer {A}bhandlung {\"u}ber die {B}ernoullischen
  {Z}ahlen.
\newblock {\em Astr. Nachr.}, 17:351--352, 1840.

\bibitem[IR90]{ireland90classical}
K.~Ireland and M.~Rosen.
\newblock {\em A {C}lassical {I}ntroduction to {M}odern {N}umber {T}heory},
  volume~84 of {\em Graduate Texts in Mathematics}.
\newblock {S}pringer-{V}erlag, 2nd edition, 1990.

\bibitem[Joh74]{Johnson74}
W.~Johnson.
\newblock Irregular prime divisors of the {Bernoulli} numbers.
\newblock {\em Mathematics of Computation}, 28(126):653--657, April 1974.

\bibitem[Kel02]{kellner02irrpairord}
B.~C. Kellner.
\newblock {\"U}ber irregul{\"a}re {P}aare h{\"o}herer \mbox{{O}rdnungen}.
\newblock {\em Diplomarbeit. Mathe\-matisches Institut der
  Georg--August--Universit{\"a}t zu G{\"o}ttingen, Germany. Also available at
  http://www.bernoulli.org/$\sim$bk/irrpairord.pdf}, 2002.

\bibitem[Kel04]{kellner04irrpowbn}
B.~C. Kellner.
\newblock On irregular prime powers of {B}ernoulli numbers.
\newblock {\em Preprint, submitted to Mathematics of Computation}, 1--29,
  2004. arXiv:math.NT/ 0409223

\bibitem[KL64]{leopoldt64zeta}
T.~Kubota and H.~W. Leopoldt.
\newblock Eine $p$-adische {T}heorie der {Z}etawerte {I}.
\newblock {\em J. Reine Angew. Math.}, 214/215:328--339, 1964.

\bibitem[Kob96]{koblitz96padic}
N.~Koblitz.
\newblock {\em $p$-adic {N}umbers, $p$-adic {A}nalysis and {Z}eta-{F}unctions},
  volume~58 of {\em Graduate Texts in Mathematics}.
\newblock {S}pringer-{V}erlag, 2nd edition, 1996.

\bibitem[Kum50]{kummer50flt}
E.~E. Kummer.
\newblock 1) {Z}wei besondere {U}ntersuchungen {\"u}ber die {C}lassen-{A}nzahl
  und {\"u}ber die {E}inheiten der aus $\lambda$-ten {W}urzeln der {E}inheit
  gebildeten complexen {Z}ahlen. 2) {A}llgemeiner {B}eweis des {F}ermat'schen
  {S}atzes, dass die {G}leichung $x^{\lambda}+y^{\lambda}=z^{\lambda}$ durch
  ganze {Z}ahlen unl{\"o}sbar ist, f{\"u}r alle diejenigen
  {P}otenz-{E}xponenten $\lambda$, welche ungerade {P}rimzahlen sind und in den
  {Z}{\"a}hlern der ersten $({\lambda}-3)/2$ {B}ernoulli'schen {Z}ahlen als
  {F}actoren nicht vorkommen.
\newblock {\em J. Reine Angew. Math.}, 40:117--129, 130--138, 1850.

\bibitem[Kum51]{kummer50bern}
E.~E. Kummer.
\newblock {\"U}ber eine allgemeine {E}igenschaft der rationalen
  {E}ntwicklungscoefficienten einer bestimmten {G}attung analytischer
  {F}unctionen.
\newblock {\em J. Reine Angew. Math.}, 41:368--372, 1851.

\bibitem[Sie64]{siegel64irr}
C.~L. Siegel.
\newblock Zu zwei {B}emerkungen {K}ummers.
\newblock {\em Nachrichten der Aka\-demie der Wissenschaften in G{\"o}ttingen.
  Mathematisch-physikalische Klasse (Gesammelte Abhandlungen, Band III,
  436--442)}, 6:51--57, 1964.

\bibitem[Tha04]{thang04adam}
R.~Thangadurai.
\newblock Adams theorem on {B}ernoulli numbers revisited.
\newblock {\em J. Number Theory}, 106:169--177, 2004.

\bibitem[vS40]{staudt40bern}
K.~G.~C. von Staudt.
\newblock Beweis eines {L}ehrsatzes die {B}ernoulli'schen {Z}ahlen betreffend.
\newblock {\em J. Reine Angew. Math.}, 21:372--374, 1840.

\bibitem[Was97]{washington97cyclo}
L.~C. Washington.
\newblock {\em Introduction to {C}yclotomic {F}ields}, volume~83 of {\em
  Graduate Texts in Mathematics}.
\newblock {S}pringer-{V}erlag, 2nd edition, 1997.

\bibitem[Yam76]{Yama76bern}
I.~Yamaguchi.
\newblock On a {B}ernoulli numbers conjecture.
\newblock {\em J. Reine Angew. Math.}, 288:168--175, 1976.

\end{thebibliography}
